\documentclass[12pt]{amsart}
\usepackage{amsmath,amscd,amssymb,amsfonts,graphics}
\newcommand{\h}{\hbox}
\newcommand{\q}{\quad}

\newcommand{\bs}{\par\bigskip}
\newcommand{\ms}{\par\medskip}
\newcommand{\sk}{\par\smallskip}
\newcommand{\bsn}{\par\bigskip\noindent}
\newcommand{\msn}{\par\medskip\noindent}

\newcommand{\1}{\hskip1pt}
\newcommand{\mcap}{\hbox{$\bigcap$}}
\newcommand{\mcup}{\hbox{$\bigcup$}}

\newcommand{\msum}{\hbox{$\sum$}}
\newcommand{\mopl}{\hbox{$\bigoplus$}}
\newcommand{\mprod}{\hbox{$\prod$}}

\newcommand{\E}{{\mathcal E}}
\newcommand{\F}{{\mathcal F}}
\newcommand{\G}{{\mathcal G}}

\newcommand{\Hc}{{\mathcal H}}
\newcommand{\Lc}{{\mathcal L}}

\newcommand{\Yc}{{\mathcal Y}}

\newcommand{\PP}{{\mathbb P}}
\newcommand{\Q}{{\mathbb Q}}
\newcommand{\C}{{\mathbb C}}

\newcommand{\RR}{{\mathbf R}}

\newcommand{\Z}{{\mathbb Z}}

\newcommand{\Gr}{{\rm Gr}}

\newcommand{\ep}{\varepsilon}
\newcommand{\ga}{\gamma}

\newcommand{\Yct}{\widetilde{\mathcal Y}}
\newcommand{\Zt}{\widetilde{Z}}
\newcommand{\zt}{\widetilde{z}}
\newcommand{\rhot}{\widetilde{\rho}}
\newcommand{\ett}{\widetilde{\eta}}
\newcommand{\lat}{\widetilde{\la}}
\newcommand{\gat}{\widetilde{\gamma}}
\newcommand{\xit}{\widetilde{\xi}}
\newcommand{\Xit}{\widetilde{\Xi}}
\newcommand{\Eo}{{}\,\overline{\!E}{}}

\newcommand{\Ff}{F_{\!f}}
\newcommand{\al}{\alpha}
\newcommand{\be}{\beta}
\newcommand{\la}{\lambda}
\newcommand{\Si}{\Sigma}
\newcommand{\si}{\sigma}

\newcommand{\lan}{\langle}
\newcommand{\ran}{\rangle}

\newcommand{\bl}{\bigl}
\newcommand{\br}{\bigr}
\newcommand{\pl}{\1{+}\1}
\newcommand{\mi}{\1{-}\1}
\newcommand{\eq}{\,{=}\,}

\newcommand{\nes}{\,{\ne}\,}
\newcommand{\stmi}{\,{\setminus}\,}
\newcommand{\less}{\,{\leqslant}\,}
\newcommand{\gess}{\,{\geqslant}\,}

\newcommand{\ssb}{\raise.15ex\h{${\scriptscriptstyle\bullet}$}}
\newcommand{\ssc}{\,\raise.15ex\h{${\scriptstyle\circ}$}\,}
\newcommand{\onto}{\twoheadrightarrow}
\newcommand{\into}{\hookrightarrow}
\newcommand{\simto}{\,\,\rlap{\hskip1.5mm\raise1.4mm\hbox{$\sim$}}\hbox{$\longrightarrow$}\,\,}
\newcommand{\plim}{\rlap{\raise-5.5pt\h{$\,\leftarrow$}}{\rm lim}}
\begin{document}
\title[Topological computation]
{Topological computation of the first Milnor fiber cohomology of hyperplane arrangements}
\author[M. Saito]{Morihiko Saito}
\address{M. Saito : RIMS Kyoto University, Kyoto 606-8502 Japan}
\email{msaito@kurims.kyoto-u.ac.jp}
\begin{abstract} We study a topological method to calculate the first Milnor fiber cohomology of a defining polynomial of a reduced projective hyperplane arrangement $X$ of degree $d$. We can show the vanishing of a monodromy eigenspace of the first Milnor fiber cohomology with eigenvalue of order $m\ge 2$ if $X\setminus(X^{[(m)]}\cup X^{\lan 3\ran})$ or more generally $X\setminus(X^{[(m)]}\cup X^{\lan 3\ran}\cup X_d)$ is connected. Here $X^{[(m)]}$ is the set of points of $X$ with multiplicity divisible by $m$, and $X^{\lan 3\ran}:=\bigcup_{i,j,k}X_i\cap X_j\cap X_k$ with $X_i$ the irreducible components of $X$, where the union is taken over $i,j,k$ with ${\rm codim}\,X_i\cap X_j\cap X_k=3$. This hypothesis can be relaxed to some extent. The assertion is reduced to the case of a line arrangement in ${\bf P}^2$ by Artin's vanishing theorem (where $X^{\lan 3\ran}=\emptyset$), and we use a projection from ${\bf P}^2$ to ${\bf P}^1$ with center a sufficiently general point of $X_d$. It may be expected that the assumption of an improved assertion is always satisfied for $m\ge 5$ (and also for $m=4$ except the Hessian arrangement). The resulting vanishing of eigenspaces has been conjectured for $m\ge 5$.
\end{abstract}
\maketitle
\ms\centerline{\bf Introduction}
\bsn
Let $f$ be a reduced polynomial defining a projective hyperplane arrangement $X\subset Y:=\PP^{n-1}$ of degree $d$. For a calculation of monodromy eigenspaces of its first Milnor fiber cohomology $H^1(\Ff,\C)_{\la}$ with eigenvalue $\la$ of order $m\gess 2$, we may assume $d/m\in\Z$, since $H^j(\Ff,\C)_{\la}\eq 0$ otherwise, where $\Ff$ denotes the Milnor fiber of $f$. These eigenspaces are calculated by the corresponding {\it Aomoto complex\1} if some hypothesis from \cite{ESV} is satisfied (although the condition is rather restrictive, and is often unsatisfied), see \cite{BDS}, \cite{BSY}, \cite{Fa}, \cite{ac}, \cite{eff}, etc.
So we consider a topological way to calculate it as in Section~1 below.
\sk
Define
\vskip-6mm
$$\aligned X^{[(m)]}&:=\{P\in X\mid{\rm mult}_PX\in m\1\Z\},\\ X^{\lan 3\ran}&:=\mcup_{i,j,k}\,X_i\cap X_j\cap X_k,\\ X'_{\lan m\ran}&:=X\stmi(X^{[(m)]}\cup X^{\lan 3\ran}\cup X_d).\endaligned$$
Here the $X_i\,\,(i\in[1,d])$ are irreducible components of $X$, and the union in the definition of $X^{\lan 3\ran}$ is taken over $i,j,k\in[1,d]$ with ${\rm codim}_YX_i\cap X_j\cap X_k\eq 3$. There is a partition
$$\h{$\bigsqcup$}_{j=1}^{r'_m}\,I^{\prime(j)}=\{1,\dots,d{-}1\},$$
such that the $X^{\prime(j)}_{\lan m\ran}:=\mcup_{i\in I^{\prime(j)}}\,X_i\cap X'_{\lan m\ran}$ ($j\in[1,r'_m]$) are connected components of $X'_{\lan m\ran}$. Set
$$I^{\prime(j)}_Q:=\{i\in[1,d{-}1]\mid X_i\ni Q,\,i\in I^{\prime(j)}\}]\q\q(Q\in X'),$$
where $X':=X\stmi X_d$. We can show the following.
\msn
{\bf Theorem~1.} {\it We have $H^1(\Ff,\C)_{\la}\eq 0$ for $\,{\rm ord}\,\la\eq m\gess 2$ if $X'_{\lan m\ran}$ is connected, or if the $X^{\prime(j)}_{\lan m\ran}$ are irreducible for $j\in[2,r'_m]$ with $r'_m\gess 2$ and ${\rm codim}_Y\,\mcap_{j=2}^{r'_m}\,\overline{X^{\prime(j)}_{\lan m\ran}}\less 2$, assuming also that $X$ is essential $($that is, it cannot be defined by a polynomial of $n{-}1$ variables$)$.}
\msn
{\bf Theorem~2.} {\it Assume for any $k\in[2,r'_m]$, there is
$$Q_k\in X^{[(m)]}\cap\overline{X^{\prime(1)}_{\lan m\ran}}\cap\overline{X^{\prime(k)}_{\lan m\ran}}\stmi(X^{\lan 3\ran}\cup X_d),$$
such that $I^{\prime(j)}_{Q_k}\eq\emptyset$ for $j\notin\{1,k\}$, and moreover $|I^{\prime(k)}_{Q_k}|$ is relatively prime to $\,m,$ and is at most $4$ in the $m$ non-prime case. Then $H^1(\Ff,\C)_{\la}\eq 0$ for $\,{\rm ord}\,\la\eq m\gess 2$.}
\ms
Theorem~1 is a special case of Theorem~2, see Remark~2.2 below. Theorems~1--2 are reduced to the case $n\eq 3$ (where $X^{\lan 3\ran}\eq\emptyset$) by Artin's vanishing theorem \cite{BBD} using an iterated general hyperplane cut. We then study the direct image of a sheaf complex calculating the vanishing cohomology under a projection to $\PP^1$ from a blow-up of $\PP^2$ along a sufficiently general point of $X_d$, see (2.2), (2.5) and Section~1 below.
\sk
A variant of Theorem~1 has been known for instance in the case $H^1(\Ff,\C)_{\la}$ is calculable by the method in \cite{ESV} mentioned above, see \cite{BDS} and also \cite{Ba}, etc. (See also (2.6) below for some explicit examples.)
We have $r'_m\eq r_m$ or $r_m\mi 1$, changing appropriately the order of the $X_i$, where $r_m$ is the number of connected components of $X_{\lan m\ran}:=X\stmi(X^{[(m)]}\cup X^{\lan 3\ran})$. This can be verified by using the {\it dual $(m)$-graph\1} of $X$. (Its vertices correspond to the irreducible components of $X$, and two vertices corresponding to $X_i,X_j$ are connected by an edge if and only if the multiplicity of $X$ at a general point of $X_i\cap X_j$ is {\it not\1} divisible by $m$.) We may have $r'_m\eq r_m\mi 1$ in the case $X_{\lan m\ran}$ has an irreducible connected component, see Remark~2.6\,(ii) below.
\sk
There is no explicit example with hypothesis of Theorem~2 unsatisfied for $m\gess 4$ (with $d/m\in\Z$) for the moment except the {\it Hessian\1} arrangement case, where the dual $(4)$-graph has 4 connected components (since it is a 4-net), see \cite{BDS}, \cite{Di2}, \cite{Yu}, and (2.6) below. By Theorems~1--2, calculations of $H^1(\Ff,\C)_{\la}$ with $m\eq{\rm ord}\,\la\gess 4$ can be reduced effectively to verifications of the hypothesis of Theorem~1 or 2 in almost all the cases. The vanishing of $H^1(\Ff,\C)_{\la}$ with ${\rm ord}\,\la\gess 5$ has been conjectured in \cite[Conjecture~1.9]{PS}, see also \cite[Problem~2]{eff}.
\sk
In Section 1 we study a topological method to describe the Milnor fiber cohomology of line arrangements. In Section 2 we prove the main theorems.
\sk
This work was partially supported by JSPS Kakenhi 15K04816.
\bs\bs
\vbox{\centerline{\bf 1. Milnor fiber cohomology of line arrangements}
\bsn
In this section we study a topological method to describe the Milnor fiber cohomology of line arrangements.}
\msn
{\bf 1.1.~Calculation of the Milnor fiber cohomology.} Let $L\subset Y:=\PP^2$ be a reduced line arrangement with $f$ a defining polynomial of degree $d$. It is well-known that there are local systems $\Lc^{k,}$ of rank 1 on $U:=Y\stmi L\,$ ($k\in[0,d{-}1]$) such that
$$H^j(\Ff,\C)_{\la}\cong H^j(U,\Lc^{k,})\q\q\bl(\la=\exp(-2\pi\sqrt{-1}\1k/d)\br),
\leqno(1.1.1)$$
and the monodromy around any irreducible component of $L$ is given by multiplication by $\la^{-1}=\exp(2\pi\sqrt{-1}\1k/d)$. This is constructed by taking the direct image of the constant sheaf under the natural morphism $f^{-1}(1)\to U$ which is an unramified covering of degree $d$, see \cite{Di1}, \cite{CS}, \cite{BS}, etc. In particular, we have $\Lc^{(0)}=\C_U$, and hence $H^j(\Ff,\C)_1=H^j(U,\C)$.
\sk
Let $\pi:\Yc\to Y:=\PP^2$ be the blow-up along a sufficiently general point $P\in L_d$ with $j_{\Yc}:U\into\Yc$ a natural inclusion. Here the $L_i$ ($i\in[1,d]$) are the irreducible components of $L$. Set
$$\F:=\RR(j_{\Yc})_*\Lc^{(k')}[2]\q\,\,\,\h{with}\q\,\,\,k'\eq d/m\,\,\,\,\,\h{or}\,\,\,\,\,d-d/m.$$
(Here we use the irreducibility of cyclotomic polynomials implicitly, see also \cite{eff}.)
This is a direct factor of the underlying $\C$-complex of a mixed Hodge module, and has the weight filtration $W$, see \cite{mhm}, \cite{BS}.
\msn
{\bf Lemma~1.1.} {\it We have the isomorphisms
$$\Gr^W_2\F={\rm IC}_{\Yc}\Lc^{(k')},\q\q\Gr^W_i\F=0\,\,\,\,(i\nes 2,3).
\leqno(1.1.2)$$
Moreover $\,\dim_{\C}\Hc^j\Gr^W_i\!\F_y=km\mi 2$ if $\,\pi(y)\in L^{[km]}$ $(k\in\Z_{>0})\,$ with $\,j\eq i\mi 3\eq{-}1\,$ or $\,0$, and it vanishes otherwise. Here $L^{[r]}:=\{P\in L\mid{\rm mult}_PL\eq r\1\}$.}
\msn
{\it Proof.} Let $\si:\Yct\to\Yc$ be the blow-up along all the points of $\pi^{-1}({\rm Sing}\,L)$ with $j_{\Yct}:U\into\Yct$ the inclusion. Then (1.1.2) holds with $\Yc$, $\F$ replaced by $\Yct$, $\RR(j_{\Yct})_*\Lc^{(k')}[2]$, and we have the isomorphism
$${\rm IC}_{\Yct}\Lc^{(k')}[-2]=(j_{\Yct})_*\Lc^{(k')}.
\leqno(1.1.3)$$
Moreover the restriction of $R^j(j_{\Yct})_*\Lc^{(k')}$ to $\si^{-1}(y)\cong\PP^1$ with $\pi(y)\in{\rm Sing}\,L$ is the direct image of a local system of rank 1 on the complement of $km$ points in $\PP^1$ (whose local monodromies are \h{{\it non-trivial}\1)} if $\pi(y)\in L^{[km]}$ with $j\eq 0$ or 1, and it vanishes otherwise. Since the Euler number of the complement of $km$ points in $\PP^1$ is $2\mi km$, we then get the last assertions of Lemma~1.1 together with the isomorphism
$$\RR\si_*{\rm IC}_{\Yct}\Lc^{(k')}={\rm IC}_{\Yc}\Lc^{(k')},
\leqno(1.1.4)$$
using the decomposition theorem \cite{BBD} (see also \cite{mhp}). This finishes the proof of Lemma~1.1.
\msn
{\bf 1.2.~Generic projection.} There is a natural projection $\rho:\Yc\to C:=\PP^1$ so that $\Yc$ is a $\PP^1$-bundle over $C$, and the exceptional divisor of the blow-up is the zero-section. (Here $C$ is identified with the set of lines on $\PP^2$ passing through $P$.) There is $c_{\infty}\in C$ such that
$$\pi\bl(\rho^{-1}(c_{\infty})\br)=L_d.$$
There are isomorphisms
$$Y':=Y\stmi L_d\cong\C^2,\q\q C':=C\stmi\{c_{\infty}\}\cong\C.$$
Let $\rho':Y'\to C'$ be the restriction of $\rho$ to $Y'\,\bl(=\Yc\stmi\pi^{-1}(L_d)\br)\cong\C^2$. Then $\rho'$ is identified with a linear function with a constant term by the above isomorphisms.
\sk
Since $P\in L_d$ is sufficiently general, we may assume the following:
$$\h{The restriction of $\rho'$ to ${\rm Sing}\,L\stmi L_d$ is {\it injective.}}
\leqno(1.2.1)$$
Set
$$C'':=C'\stmi\rho'({\rm Sing}\,L\stmi L_d).$$
\msn
{\bf Lemma~1.2.} {\it We have the vanishing
$$R^j\rho_*\Gr^W_i\!\F=0\,\,\,\,\,\h{unless $\,(i,j)\eq(2,-1)\,$ or $\,(3,0)$}.
\leqno(1.2.2)$$
For $c\in C'\cong\C$, there are equalities}
$$\dim_{\C}(R^{-1}\rho_*\Gr^W_2\!\F)_c=\begin{cases}d\mi2&\h{\it if}\,\,\,\,c\in C''\\ d\mi 3&\h{\it if}\,\,\,\,c\in\rho'(L^{[k]}),\,\,k/m\in\Z,\\ d\mi k\mi 1&\h{\it if}\,\,\,\,c\in\rho'(L^{[k]}),\,\,k/m\notin\Z,\end{cases}
\leqno(1.2.3)$$
\msn
{\it Proof.} The assertion for $\Gr^W_3\!\F$ follows from Lemma~1.1. For $\Gr^W_2\!\F$, it is enough to calculate
$$H^j(\Yc_c,{\rm IC}_{\Yc}\Lc^{(k')}|_{\Yc_c})=H^j(\Yct_c,{\rm IC}_{\Yct}\Lc^{(k')}|_{\Yct_c}),
\leqno(1.2.4)$$
using (1.1.4) (since $\rho$, $\si$ are proper), where $\Yc_c:=\rho^{-1}(c)$, $\Yct_c:=\si^{-1}(\Yc_c)$. We see that (1.2.4) vanishes for $j\nes 1$, and its dimension is as in (1.2.3) by an argument similar to the proof of Lemma~1.1. This finishes the proof of Lemma~1.2.
\msn
{\bf 1.3.~Description of the first Milnor fiber cohomology.} Let $j_{C'}:C'\into C$ be a natural inclusion. Set
$$\G:=\RR\rho_*\F=\RR(j_{C'})_*\G'\q\h{with}\q\G':=\G|_{C'}.$$
(This is equal to $^{\bf p}\Hc^0\RR\rho_*\F\,$ in the notation of \cite{BBD}.) We have the weight filtration $W$ on $\G'$ with
$$\Gr^W_i\G'=(\RR\rho_*\Gr^W_i\!\F)|_{C'}\q\q(i\eq2,3).
\leqno(1.3.1)$$
We have
$${\rm Supp}\,\Gr^W_3\G'\subset C'\stmi C'',\q\Hc^j\Gr^W_3\G'=0\,\,\,\,(j\nes 0),
\leqno(1.3.2)$$
and $\Gr^W_2\G'$ is an intersection complex with coefficients in a local system. More precisely,
$\Gr^W_2\G'[-1]$ is a sheaf in the usual sense on $C'$, and is the open direct image of a local system $\G''$ of rank $d{-}2$ on $C''$, where $\G'':=\G'[-1]|_{C''}$. (This can be shown using the decomposition theorem \cite{BBD} or more directly by Lemma~1.4 below.)
\sk
We thus get the following.
\msn
{\bf Proposition~1.3.} {\it There are isomorphisms
$$(R^{-1}\rho_*\Gr^W_2\!\F)|_{C'}=\Hc^{-1}\G',
\leqno(1.3.3)$$
\vskip-6mm
$$H^1(\Ff,\C)_{\la}\cong H^{-1}(C',\G')=\Gamma(C'',\G''),
\leqno(1.3.4)$$
where $\la=\exp(\pm 2\pi\sqrt{-1}/m)$. $($Note that the monodromy is defined over $\Q.)$}
\msn
{\bf 1.4.~Description of the vanishing cycles.} We have the following.
\msn
{\bf Lemma~1.4.} {\it Let $c\in\rho'(L^{[k]})\stmi\{c_{\infty}\}\,\subset\,C'\cong\C$. Then
$$\dim_{\C}\varphi_{t-c,\,1}\,\G''=1,\,\,\,\,\dim_{\C}\varphi_{t-c,\,\ne 1}\,\G''=0\q\h{if}\q k/m\in\Z.
\leqno(1.4.1)$$
\vskip-7mm
$$\dim_{\C}\varphi_{t-c,\,1}\,\G''=0,\,\,\,\,\dim_{\C}\varphi_{t-c,\,\ne 1}\,\G''=k{-}1\q\h{if}\q k/m\notin\Z.
\leqno(1.4.2)$$
Here $\,\varphi_{t-c,\,1}$ and $\,\varphi_{t-c,\,\ne 1}\,$ denote respectively the unipotent and non-unipotent monodromy part of the vanishing cycle functor $\varphi_{t-c}$ with $t$ the coordinate of $C'\cong\C$.}
\msn
{\it Proof.} This is shown by calculating the vanishing cycles $\varphi_{\rhot^*t-c}(j_{\Yct})_*\Lc^{(k')}$ on $\Yct$ using (1.1.3--4) together with the commutativity of the nearby and vanishing cycle functors with the direct image under a proper morphism, where $\rhot:=\rho\ssc\si$.
\sk
In the case $k/m\in\Z$, the direct image sheaf $(j_{\Yct})_*\Lc^{(k')}$ is a local system on the complement of the proper transform of $L$ in an open neighborhood of $\,\si^{-1}(y)$ with $\pi(y)\in L^{[k]}$. So (1.4.1) follows.
\sk
In the case $k/m\notin\Z$, the restriction of $\varphi_{\rhot^*t-c,\,\ne 1}(j_{\Yct})_*\Lc^{(k')}$ to $\,\si^{-1}(y)$ with $\pi(y)\in L^{[k]}$ is the direct image of a local system of rank 1 on the complements of $k{+}1$ points in $\PP^1$ (where the local monodromies are non-trivial), and we have the vanishing
$$\varphi_{\rhot^*t-c,\,1}(j_{\Yct})_*\Lc^{(k')}=0,$$
since we can verify the isomorphism
$$(j_{\Yct})_*\Lc^{(k')}|_{\{\rhot^*t=c\}}\simto\psi_{\rhot^*t-c,\,1}(j_{\Yct})_*\Lc^{(k')}.$$
So we get (1.4.2). This finishes the proof of Lemma~1.4.
\msn
{\bf 1.5.~Description of a general stalk of $\G''$.} Let $c\in C''$. Set
$$Z_c:=\pi\bl(\rho^{-1}(c)\br)\subset Y,\q Z'_c:=Z_c\stmi L,\q\E'_c:=\Lc^{(k')}|_{Z'_c},$$
so that
$$\G''_{c}=H^1(Z'_c,\E'_c).
\leqno(1.5.1)$$
Since $Z_c\cong\PP^1$ and $d/m\in\Z$, we have a ramified covering
$$p_c:\Zt_c\to Z_c$$
of degree $m$, which is unramified over $Z'_c$, and is {\it totally ramified\1} over
$$\Xi_c:=Z_c\cap L,$$
that is, $|p_c^{-1}(z)|=1$ if $z\in\Xi_c$. Here $\Zt_c$ is normal, and hence smooth. Let
$$p'_c:\Zt'_c\to Z'_c$$
be the restriction of $p_c$ over $Z'_c$, which is an unramified covering of degree $m$. We then get an isomorphism
$$\E'_c\cong\bl((p'_c)_*\C_{\Zt'_c}\br)^{\ga^*=\la'},
\leqno(1.5.2)$$
with right-hand side the direct factor of $(p'_c)_*\C_{\Zt'_c}$ on which the action of a generator $\ga$ of the covering transformation group of $p_c$ coincides with multiplication by $\la':=\exp(2\pi\sqrt{-1}\1/m)$ (replacing $\ga$ if necessary). Here $\ga^*$ is semi-simple, since it has finite order. (Note that a local system of rank 1 on a Zariski-open subset of $\PP^1$ is determined by the local monodromies.)
\sk
Since $p_c$ is totally ramified at $\Xi_c$, we get by (1.5.1--2) the following.
\msn
{\bf Proposition~1.5.} {\it We have the isomorphism}
$$\G''_c\cong H^1(\Zt_c,\C)^{\ga^*=\la'}.
\leqno(1.5.3)$$
\sk
(The right-hand side denotes the direct factor of $H^1(\Zt_c,\C)$ on which the action of $\ga^*$ coincides with multiplication by $\la'$. Note that the isomorphisms in (1.5.2--3) are canonical up to constant multiplication by $\C^*$.)
\bs\bs
\vbox{\centerline{\bf 2. Proof of Theorems~1 and 2}
\bsn
In this section we prove the main theorems.}
\msn
{\bf 2.1.~Description of the right-hand side of {\rm(1.5.3)}.} We can describe the right-hand side of (1.5.3) using the isomorphism by Poincar\'e duality:
$$H^1(\Zt_c,\C)=H^1(\Zt_c,\C)^{\vee}=H_1(\Zt_c,\C).
\leqno(2.1.1)$$
Note that for $\tau,\tau'\in H^1(\Zt_c,\C)$, $\si\in H_1(\Zt_c,\C)$, we have
$$\lan\ga^*\tau,\ga^*\tau'\ran=\lan\tau,\tau'\ran,\q\lan\ga^*\tau,\si\ran=\lan\tau,\ga_*\si\ran.
\leqno(2.1.2)$$
\sk
Let $z_i\in Z_c$ with $\{z_i\}=Z_c\cap L_i$ for each line $L_i\subset L$ (where $z_d=P$). Note that $\Xi_c=\{z_1,\dots,z_d\}$. Choose a path $\xi_i$ between $z_d$ and $z_i$, which is contained in $Z'_c$ except the both ends. Here we may assume $\xi_i\cap\xi_j=\{z_d\}$ ($i\nes j$). We have the following.
\msn
{\bf Lemma~2.1.} {\it Choosing $z_{\infty}\in Z_c$, its complement $Z_c\stmi\{z_{\infty}\}$ is identified with $\C$, and the $\xi_i$ can be defined as the closed line segment with endpoints $z_i$, $z_d$ in the $2$-dimensional real affine space $\C$, if $P,c$ are sufficiently general.}
\msn
{\it Proof.} Let $\ell_i$ be the half line in the real 2-dimensional affine space $Z_c\stmi\{z_{\infty}\}\cong\C$ with endpoint $z_d$ and containing $z_i$. (Recall that any algebraic automorphism of the complex algebraic variety $\C\subset\PP^1$ is given by $\al\1x\pl\be$ with $\al\in\C^*,\be\in\C$, where $x$ is an algebraic coordinate of $\C\subset\PP^1$. This says also that the ambiguity of algebraic coordinate of $\C\subset\PP^1$ is given by it.) We have to show that $\ell_i\nes\ell_j$ ($i\nes j$) if $P,c$ are sufficiently general.
\sk
Take coordinates $x,y$ of $\C^2\cong\PP^2\stmi L_d$ such that the equation of $L'_i:=L_i\stmi L_d\subset\C^2$ ($i\in[1,d{-}1]$) can be written as
$$y=a_ix-b_i\q(a_i,b_i\in\C),$$
that is, $L'_i\cap\{x\eq 0\}\nes\emptyset$. We consider only $P,c$ such that the equation of $Z_c\stmi\{z_d\}=\rho'{}^{-1}(c)$ is written as
$$y=\al\1x-\be\q(\al,\be\in\C).$$
Here the ambiguity (or freedom) of $P,c$ corresponds to that of $\al,\be$ respectively. We will show that $\ell_i\nes\ell_j$ ($i\nes j$) if $\be$ is sufficiently general for each $\al$ with (1.2.1) satisfied.
\sk
We define $z_{\infty}$ to be the intersection of $Z_c\stmi\{z_d\}$ with $\{x\eq 0\}\subset\C^2$. Then $x^{-1}$ can be extended to a complex coordinate of $Z_c\stmi\{z_{\infty}\}\cong\C$ with origin $z_d$. So it is enough to show that
$$\arg x(z_i)\nes\arg x(z_j)\q(i\nes j),
\leqno(2.1.3)$$
where $x(z_i)$ is the value of the coordinate $x$ at $z_i$. Put
$$g_i(\be):=\arg(\be\mi b_i),\q h_i(\al):=\arg(\al\mi a_i).$$
Since $\{z_i\}=L_i\cap Z_c$, we get
$$x(z_i)=(\be\mi b_i)/(\al\mi a_i)$$
from the above equations. So condition (2.1.3) is equivalent to that
$$g_i(\be)\mi h_i(\al)\ne g_j(\be)\mi h_j(\al)\q(i\nes j),
\leqno(2.1.4)$$
assuming that $\be$ is sufficiently general. Here ``sufficiently general" means that the set of $\be$ satisfying (2.1.4) is a dense open subset of $\C$. (Recall that $\,\arg x\,$ for $x\in\C$ is the imaginary part of a brach of $\log x$ contained in $[0,2\pi)$.)
\sk
Set $\,\Si_{i,r_i}:=\{g_i(\be)\eq r_i\}\subset\C$ for $r_i\in[0,2\pi)$. (Note that $\Si_{i,r_i}$ is a real half line with endpoint $b_i$.) Let $U_{\be}\subset\C$ be the complement of the union of real lines $\ell_{i,j}$ in $\C$ passing through $b_i,b_j$ ($i\nes j$). Then $\Si_{i,r_i}\cap\Si_{j,r_j}\cap U_{\be}$ is discrete, and $\Si_{j,r_j}\cap U_{\be}$ is not contained in $\Si_{i,r_i}\cap U_{\be}$ for $r_i,r_j\in[0,2\pi)$. This implies that we can change the left-hand side of (2.1.4) slightly by changing $\be$ with right-hand side unchanged. We can then get (2.1.4) by induction of the number of the $i$ such that (2.1.4) does not hold for some $j$. Each change of $\be$ becomes much smaller and smaller as induction goes up so that once we get (2.1.4) for some $i,j$, we never loose it by changing $\be$ further. This finishes the proof of Lemma~2.1.
\ms
Using $\xi_i$ in Lemma~2.1 together with a sufficiently small circle around $z_i$, we get an element
$$\eta_i\in \pi_1(Z_c\stmi\Xi'_c,z_d)\q\h{with}\q\Xi'_c:=\Xi_c\stmi\{z_d\}.
\leqno(2.1.5)$$
This can be lifted non-canonically to an element
$$\ett_i\in \pi_1(\Zt_c\stmi\Xit'_c,\zt_d)\q\h{with}\q p_c(\Xit'_c)=\Xi_c,\,\,p_c(\zt_d)=z_d.
\leqno(2.1.6)$$
Its ambiguity is given by the action of the covering transformation group $\Gamma_{p_c}\cong\mu_m$ of $p_c$. Set $\la'':=\la'{}^{-1}$. We denote by $\si_i$ the image of $\ett_i$ in $H_1(\Zt_c,\C)$ (viewed as a topological cycle), and by $\si^{\la''}_i$ its image in $H_1(\Zt_c,\C)^{\ga_*=\la''}$ using the eigenvalue decomposition by the action of $\ga_*$ on $H_1(\Zt_c,\C)$. The ambiguity of $\si^{\la''}_i$ is then given by the natural action of $\mu_m\subset\C^*$ (by scalar multiplication), since the isomorphism $\Gamma_{p_c}\cong\mu_m$ is given by the generator $\ga\in\Gamma_{p_c}$. We have the following.
\msn
{\bf Proposition~2.1.} {\it For $\,j\in[1,d{-}1]$, set $\,I_{(j)}:=\{1,\dots,d{-1}\}\stmi\{j\}$. Then the $\,\si^{\la''}_i\,\,(i\in I_{(j)})$ form a $\C$-basis of $H_1(\Zt_c,\C)^{\ga_*=\la''}\,(=H^1(\Zt_c,\C)^{\ga^*=\la'})$.}
\msn
{\it Proof.} Since $\dim H_1(\Zt_c,\C)^{\ga_*=\la''}\eq d\mi 2$ by (1.5.3), the assertion is reduced to showing that the $\si^{\la''}_i$ ($i\nes j$) span $H_1(\Zt_c,\C)^{\ga_*=\la''}$ for any $j\in[1,d{-}1]$. For this it is enough to prove that the $\si_i$ ($i\nes j$) generate $H_1(\Zt_c,\C)$ over the group ring $\C[\Gamma_{p_c}]$. (Note that $\Gamma_{p_c}\cong\mu_m$ is generated by $\ga$.) Here we have the surjection
$$H_1(\Zt_c\stmi\{\zt_j\},\C)\onto H_1(\Zt_c,\C).
\leqno(2.1.7)$$
We see that $\bigcup_{i\ne j}\xi_i$, $\bigcup_{i\ne j}\xit_i$ are deformation retracts of $Z_c\stmi\{z_j\}$ and $\Zt_c\stmi\{\zt_j\}$ respectively, where $\xit_i:=p_c^{-1}(\xi_i)$.
(Indeed, choosing $z_{\infty}\in Z_c$, its complement $Z_c\stmi\{z_{\infty}\}$ is identified with $\C$, and the $\xi_i$ can be obtained as the convex hull of $z_i$ and $z_d$ in the real vector space $\C$ by Lemma~2.1. Consider the real half line $\ell_j$ in $\C$ which contains $z_j$ as the endpoint and such that the real full line containing $\ell_j$ contains $z_d$ although $z_d\notin\ell_j$. Then the complement in $Z_c$ of an appropriate open neighborhood of the closure of $\ell_j$ in $Z_c$ is a deformation retract of $Z_c\stmi\{z_j\}$, and $\bigcup_{i\ne j}\xi_i$ is a deformation retract of the complement. The assertion for $\bigcup_{i\ne j}\xit_i$, $\Zt_c\stmi\{\zt_j\}$ follows from this, since $\Zt_c$ is an unramified covering over $Z'_c$.)
The assertion then follows, since $\mu_m$ acts as the covering transformation group of the unramified covering
$$\xit_i\stmi\{\zt_i,\zt_d\}\to\xi_i\stmi\{z_i,z_d\}.$$This finishes the proof of Proposition~2.1.
\msn
{\bf Corollary~2.1.} {\it The $\si^{\la''}_i\,\,(i\in[1,d{-}1])$ span $H_1(\Zt_c,\C)^{\ga_*=\la''}$, and there are $\al_i\in\C^*$ with
$$\msum_{i=1}^{d-1}\,\al_i\1\si^{\la''}_i=0.
\leqno(2.1.8)$$
Hence we get the isomorphism
$$H_1(\Zt_c,\C)^{\ga_*=\la''}=\bl(\mopl_{i=1}^{d-1}\,\C\1e_i\br)/\C\bl(\msum_{i=1}^{d-1}\,\al_i\1e_i\br),
\leqno(2.1.9)$$
where $(e_1,\dots,e_{d-1})$ is the basis of the complex vector space $\C^{d-1}$, and $\si^{\la''}_i$ is identified with the image of $e_i$ in the right-hand side of $(2.1.9)$ for $i\in\{1,\dots,d{-}1\}$.}
\msn
{\bf 2.2.~Proof of Theorem~1 in the $X'_{\lan m\ran}$ connected case.} The assertion is reduced to the line arrangement case by Artin's vanishing theorem \cite{BBD} using an iterated general hyperplane cut. So $X$ will be denoted by $L$ and the notation in Section~1 will be used. By Proposition~1.3, it is enough to show that
$$\Gamma(C'',\G'')=0,\,\,\,\h{or equivalently,}\,\,\,(\G''_c)^G=0.
\leqno(2.2.1)$$
Here $c\in C''$ as in (1.5), and $G$ is the monodromy group of the local system $\G''$ on $C''$ with base point $c$.
\sk
By Proposition~1.5 and Corollary~2.1 together with Poincar\'e duality as in (2.1.1--2), any element of $(\G''_c)^G$ can be expressed by
$$\msum_{i=1}^{d-1}\,\be_i\1\si^{\la''}_i,$$
where the $\be_i\in\C$ are unique up to addition by $\ep\1\al_i$ with $\ep\in\C$. Assuming $X'_{\lan m\ran}$ is connected, we have to show that
$$\be_i/\al_i\,\,\,\h{is independent of}\,\,\,i\in\{1,\dots,d{-}1\}.
\leqno(2.2.2)$$
\sk
Here $\{1,\dots,d{-}1\}$ is identified with the vertices of the dual $(m)$-graph. Since the latter is connected by hypothesis, it is enough to show (2.2.2) in the case $i$ belong to a subset $I_z$ such that the $L_i$ ($i\in I_z$) are the lines in $L$ passing through some $z\in L^{[k]}\stmi L_d$ with $k/m\notin\Z$, where $|I_z|\eq k$.
\sk
Choosing a path between $c$ and $\rho'(z)$, and taking a sufficiently small circle around $\rho'(z)$, we get $\zeta\in\pi_1(C'',c)$ going around $\rho'(z)$. Consider the short exact sequence
$$0\to(\G''_c)^{\zeta}\buildrel{\iota}\over\to\G''_c\to{\rm Coker}\,\iota\to 0,
\leqno(2.2.3)$$
where the first term is the invariant part of $\G''_c$ under the monodromy action by $\zeta$.
Using (1.2.3) with (1.3.3) and (1.4.2), the dimensions of the three terms are given respectively by
$$d\mi k\mi 1,\q d\mi 2,\q k\mi 1.$$
We can verify that the $\si^{\la''}_i$ for $i\in I_z^c:=\{1,\dots,d{-}1\}\stmi I_z\,$ belong to $(\G''_c)^{\zeta}$. Then they span it using Proposition~2.1, since $|I_z^c|\eq d\mi k\mi 1$. We thus get by Corollary~2.1 the isomorphism
$${\rm Coker}\,\iota=\bl(\mopl_{i\in I_z}\,\C\1e_i\br)/\C\bl(\msum_{i\in I_z}\,\al_i\1e_i\br),
\leqno(2.2.4)$$
such that the image of $\si^{\la''}_i$ in ${\rm Coker}\,\iota$ is identified with the image of $e_i$ in the right-hand side of (2.2.4) for $i\in I_z$.
\sk
Since $\msum_{i=1}^{d-1}\,\be_i\1\si^{\la''}_i\in(\G''_c)^G$ belongs to $(\G''_c)^{\zeta}$, its image in ${\rm Coker}\,\iota$ vanishes. This means that
$$\msum_{i\in I_z}\,\be_i\1e_i\in\C\bl(\msum_{i\in I_z}\,\al_i\1e_i\br).$$
So we get (2.2.2) for $i\in I_z$. This finishes the proof of Theorem~1 in the $X'_{\lan m\ran}$ connected case.
\msn
{\bf Remark~2.2.} Theorem~1 is a special case of Theorem~2. Indeed, this can be verified easily for $r'_m\less 2$, since $X$ is essential. In the case $r'_m\gess 3$, the latter property of $X$ implies that there is a hyperplane $X_i\subset\overline{X^{\prime(1)}_{\lan m\ran}}$ not containing the projective subspace $\mcap_{j=2}^{r'_m}\,\overline{X^{\prime(j)}_{\lan m\ran}}$ of codimension~2 so that $X_i\cap\overline{X^{\prime(j)}_{\lan m\ran}}$ does not coincide with it for each $j\in[2,r'_m]$. The assumption of Theorem~2 is then satisfied.
\msn
{\bf 2.3.~Construction of the covering space $\Zt_c$.} Let $\xi_i$ be as in Lemma~2.1. We may assume
$$\h{The $\arg z(z_i)$ are strictly increasing ($i\in[1,d{-}1]$).}
\leqno(2.3.1)$$
Here $z$ is a complex coordinate of $Z_c\stmi\{z_{\infty}\}\cong\C$ with $z(z_d)\eq 0$. Set
$$E_k:=Z_c\stmi\mcup_{i=1}^{d-1}\,\xi_i\q(k\in\Z/m\1\Z).$$
There is a compactification $\Eo_k$ of $E_k$ together with a surjection $\phi_k:\Eo_k\to Z_c$ such that the boundary is as follows:
$$\Eo_k\stmi E_k=\mcup_{i=1}^{d-1}\,\xi^{k,+}_i\,\cup\,\mcup_{i=1}^{d-1}\,\xi^{k,-}_i,$$
where $\xi^{k,+}_i$, $\xi^{k,-}_i$ are isomorphic to $\xi_i$ by $\phi_k$. Moreover
$$\xi^{k,+}_i\cap\xi^{k,-}_i=\{p^k_i\},\q\xi^{k,+}_i\cap\xi^{k,-}_{i+1}=\{q^k_i\},$$
with $\phi_k(p^k_i)\eq z_i$, $\phi_k(q^k_i)\eq z_d$ for $i\in\{1,d{-}1\}\cong\Z/(d{-}1)\Z$. The picture of the images of $\xi^{k,\pm}_i$, $p^k_i$, etc.\ by $\phi_k$ is as follows (where $\phi_k$ is omitted to simplify the picture):
$$\setlength{\unitlength}{5mm}
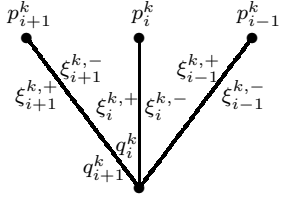
\begin{picture}(7,5.2)
\linethickness{0.3mm}
\qbezier(4,0)(2.5,2)(1,4)
\qbezier(4,0)(4,2)(4,4)
\qbezier(4,0)(5.5,2)(7,4)
\put(4,4){\circle*{.3}}
\put(1,4){\circle*{.3}}
\put(7,4){\circle*{.3}}
\put(4,0){\circle*{.3}}
\put(.5,4.5){$\scriptstyle p^k_{i+1}$}
\put(3.8,4.5){$\scriptstyle p^k_i$}
\put(6.6,4.5){$\scriptstyle p^k_{i-1}$}
\put(4.15,2){$\scriptstyle \xi^{k,-}_i$}
\put(2.85,2){$\scriptstyle \xi^{k,+}_i$}
\put(.7,2.3){$\scriptstyle \xi^{k,+}_{i+1}$}
\put(1.9,3){$\scriptstyle \xi^{k,-}_{i+1}$}
\put(5,3){$\scriptstyle \xi^{k,+}_{i-1}$}
\put(6.2,2.3){$\scriptstyle \xi^{k,-}_{i-1}$}
\put(3.37,.98){$\scriptstyle q^k_i$}
\put(2.5,.4){$\scriptstyle q^k_{i+1}$}
\end{picture}$$
\sk
We get the covering space $\Zt_c$ topologically by identifying $\,\xi^{k,+}_i$ with $\,\xi^{k+1,-}_i$; in particular, $p^k_i$ with $p^k_i$, and $q^k_i$ with $q^{k+1}_{i-1}$. (Related to the latter, note that $m$ and $d{-}1$ are mutually prime, hence $(1,-1)\in(\Z/m\1\Z){\times}(\Z/(d{-}1)\1\Z)$ is a generator of the finite cyclic group.)
\msn
{\bf Remark~2.3.} By the above construction, we can choose the $\ett_i$ in (2.1.6) appropriately so that
$$\msum_{i=1}^{d-1}\,\si_i=0\q\h{in}\,\,\,H_1(\Zt_c,\C),
\leqno(2.3.2)$$
since the $E_k$ are contractible. In (2.1.8) we may then assume
$$\al_i\eq 1\q(\forall\,i\in[1,d{-}1]).
\leqno(2.3.3)$$
\msn
{\bf 2.4.~Calculation of local monodromies of $\G''$.} Let $Q\in{\rm Sing}\,L\setminus L_d$. There are coordinates $x',y'$ of $Y\setminus L_d\cong\C^2$ with origin $Q$ such that $x'\eq x\mi a$ with $x$ as in the proof of Lemma~2.1 and $y'$ is the pull-back of a coordinate of $C'\cong\C$ by $\rho'$ in (1.2). Then the $L_i\subset L$ ($i\in I_Q$) passing through $Q$ can be described as
$$x'=a'_i\1y'\q(a'_i\in\C,\,\,i\in I_Q).
\leqno(2.4.1)$$
\sk
We assume $c\in C''$ is sufficiently near $\rho'(Q)$ so that
$$\max(I_Q)-\min(I_Q)=|I_Q|\mi 1,
\leqno(2.4.2)$$
that is, $I_Q=\Z\cap[\min(I_Q),\max(I_Q)]$. (Note that (2.3.1) is assumed.)
\sk
Substituting $y'\eq\ep e^{i\theta}$ in (2.4.1) with $0\,{<}\,\ep\,{\ll}\,1$ fixed and $\,\theta\in[0,2\pi]$ varying, we can verify that, under the local monodromy around $\rho'(Q)\in C'$, the loop $\,\eta_i\in\pi_1(\Zt_c\stmi\Xi'_c,z_d)\,$ in (2.1.5) is changed as below (where we assume $\,i\eq 2\,$ with $I_Q=\{1,2,3\}\,$ for simplicity):
$$\h{$\setlength{\unitlength}{5mm}
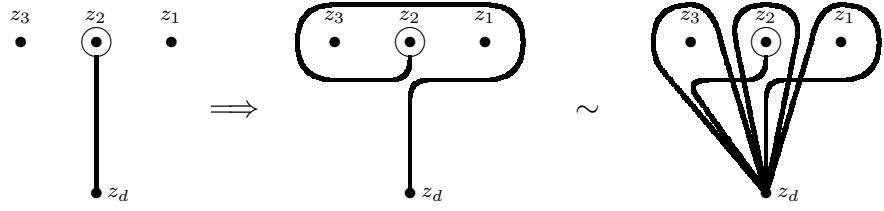
\begin{picture}(7,5.6)
\linethickness{.5mm}
\qbezier(4,0)(4,2)(4,3.6)
\put(4,4){\circle{.8}}
\put(4,4){\circle*{.3}}
\put(2,4){\circle*{.3}}
\put(6,4){\circle*{.3}}
\put(4,0){\circle*{.3}}
\put(1.7,4.6){$\scriptstyle z_3$}
\put(3.7,4.6){$\scriptstyle z_2$}
\put(5.7,4.6){$\scriptstyle z_1$}
\put(4.3,-.1){$\scriptstyle z_d$}
\end{picture}$}\raise1cm\h{$\Longrightarrow$}
\h{$\setlength{\unitlength}{5mm}
\begin{picture}(7,5.4)
\linethickness{.5mm}
\qbezier(4,0)(4,1)(4,2.5)
\qbezier(4,2.5)(4,3)(4.5,3)
\qbezier(6,3)(5,3)(4.5,3)
\qbezier(6,3)(7,3)(7,4)
\qbezier(6,5)(7,5)(7,4)
\qbezier(6,5)(4,5)(2,5)
\qbezier(2,3)(1,3)(1,4)
\qbezier(2,5)(1,5)(1,4)
\qbezier(2,3)(3,3)(3.5,3)
\qbezier(4,3.6)(4,3)(3.5,3)
\put(4,4){\circle{.8}}
\put(4,4){\circle*{.3}}
\put(2,4){\circle*{.3}}
\put(6,4){\circle*{.3}}
\put(4,0){\circle*{.3}}
\put(1.7,4.6){$\scriptstyle z_3$}
\put(3.7,4.6){$\scriptstyle z_2$}
\put(5.7,4.6){$\scriptstyle z_1$}
\put(4.3,-.1){$\scriptstyle z_d$}
\end{picture}$}\q\,\,\,\,\raise1cm\h{$\sim$}\,\,\,
\h{$\setlength{\unitlength}{5mm}
\begin{picture}(7,5.4)
\linethickness{.5mm}
\qbezier(4,0)(4,1)(4,2.5)
\qbezier(4,2.5)(4,3)(4.5,3)
\qbezier(6,3)(5,3)(4.5,3)
\qbezier(6,3)(7,3)(7,4)
\qbezier(6,5)(7,5)(7,4)
\qbezier(6,5)(5.4,5)(5.2,4)
\qbezier(4,0)(4,0)(5.2,4)
\qbezier(4,5)(5,5)(4.8,4)
\qbezier(4,0)(4,0)(4.8,4)
\qbezier(4,5)(3,5)(3.2,4)
\qbezier(4,0)(4,0)(3.2,4)
\qbezier(2,5)(2.6,5)(2.8,4)
\qbezier(4,0)(4,0)(2.8,4)
\qbezier(2,5)(1,5)(1,4)
\qbezier(1.2,3.3)(1,3.6)(1,4)
\qbezier(1.2,3.3)(4,0)(4,0)
\qbezier(2.2,2.5)(4,0)(4,0)
\qbezier(2.2,2.5)(1.8,3)(2.5,3)
\qbezier(2.5,3)(3,3)(3.5,3)
\qbezier(4,3.6)(4,3)(3.5,3)
\put(4,4){\circle{.8}}
\put(4,4){\circle*{.3}}
\put(2,4){\circle*{.3}}
\put(6,4){\circle*{.3}}
\put(4,0){\circle*{.3}}
\put(1.7,4.6){$\scriptstyle z_3$}
\put(3.7,4.6){$\scriptstyle z_2$}
\put(5.8,4.6){$\scriptstyle z_1$}
\put(4.3,-.1){$\scriptstyle z_d$}
\end{picture}$}$$
\msn
Here a thick path means a ``double path" consisting of two paths with opposite directions. These have two different liftings to the covering space, since they are connected to a circle around $z_2$.
\sk
Crossing the images of $\xi_i^{k,+}$, $\xi_i^{k,-}$ from left to right corresponds to the action of a generator $\gat$ of the covering transformation group. We thus get the following.
\msn
{\bf Lemma~2.4.} {\it Let $i_0:=\min I_Q$, $q:=|I_Q|$. The image of $\si_{i_0+i}\,\,(i\in[0,k-1])$ by the local monodromy around $\rho'(Q)$ is given by}
$$\gat_*^q\si_{i_0+i}+\gat_*^i(1\mi\gat_*)\bl(\msum_{j=0}^{q-1}\,\si_{i_0+j}\br).
\leqno(2.4.3)$$
\msn
{\bf 2.5.~Proof of Theorem~2.} Assume $r'_m\gess 2$. For $k\in[2,r'_m]$, let $Q_k$, $I^{\prime(j)}_{Q_k}$ be as in Theorem~1. By assumption, we have
$$I^{\prime(j)}_{Q_k}\eq\emptyset\,\,\,\,\,\h{if}\,\,\,\,\,j\nes 1,k,
\leqno(2.5.1)$$
\vskip-7mm
$$\h{$|I^{\prime(k)}_{Q_k}|$ is relatively prime to $m$, and is at most 4 in the $m$ non-prime case.}
\leqno(2.5.2)$$
\sk
Let $c\in C''$ be sufficiently near $\rho'(Q_k)$ as in (2.4). Construct the $\si_i$ as in (2.3) so that (2.3.3) holds.
Assume $\msum_{i=1}^{d-1}\be_i\si_i^{\la''}\in(\G''_c)^G$ as in (2.2). By the same argument as there, there are $\be^{(j)}\in\C$ ($j\in[1,r'_m]$) such that
$$\be_i=\be^{(j)}\q\h{if}\q i\in I^{\prime(j)}\q(j\in[1,r'_m]),
\leqno(2.5.3)$$
where the $I^{\prime(j)}$ are as in the introduction. Note that $\al_i\eq 1$ ($\forall\,i)$ by (2.3.3).
\sk
Replacing $\be_i$ with $\be_i\mi\be^{(1)}$ ($\forall\,i\in[1,d{-}1]$), we may assume
$$\be_i\eq\be^{(1)}\eq 0\q\h{if}\q i\in I^{\prime(1)}.
\leqno(2.5.4)$$
Using Lemma~2.4 and (2.5.1), the image of $\msum_{i=1}^{d-1}\be_i\si_i^{\la''}$ under the local monodromy around $\rho'(Q_k)$ is then given by
$$\msum_{i=1}^{d-1}\be_i\si_i^{\la''}+\be^{(k)}\bl(\msum_{i\in I^{\prime(k)}_{Q_k}}\,\lat^i\br)(1\mi\lat)\bl(\msum_{j=0}^{q-1}\,\si_{i_0+j}^{\la''}\br),
\leqno(2.5.5)$$
where $\lat$ is a primitive $m$\1th root of unity. Note that the $\si_i$ ($i\notin I_{Q_k}$) are invariant by the local monodromy around $\rho'(Q_k)$.
\sk
Assume $\be^{(k)}\nes 0$. Since $\msum_{i=1}^{d-1}\be_i\si_i^{\la''}$ is monodromy invariant by assumption, we can then deduce the following from (2.5.5) using Proposition~2.1:
$$\msum_{i\in I^{\prime(k)}_{Q_k}}\,\lat^i=0.
\leqno(2.5.6)$$
However, this contradicts (2.5.2) by Remark~2.5 below. This finishes the proof of Theorem~2.
\msn
{\bf Remark~2.5.} Let $\la$ be a primitive $m$\1th root of unity. Let $J$ be a subset of $\{0,\dots,m{-}1\}$ such that $|J|$ is relatively prime to $m$, and is at most 4 in the $m$ non-prime case. Then we have $\msum_{j\in J}\,\la^j\nes 0$. This does not hold without the last hypothesis, for instance in the case $J=\{0,3,4,8,9\}$ with $m\eq 12$, $|J|\eq 5$. (In the $m$ prime case, $\msum_{i=0}^{m-1}\,x^i\in\Q[x]$ is an irreducible polynomial by the theory of cyclotomic polynomials.)
\msn
{\bf 2.6.~Generalized Hessian arrangements.} For an integer $b\gess 2$, set $L:=\{f\eq0\}\subset\PP^2$ with
\vskip-8mm
$$f=xyz\,\mprod_{i,j=0}^{b-1}\,(\zeta^ix\pl\zeta^jy\pl z).$$
with $\zeta:=e^{2\pi\sqrt{-1}/b}$. Here ${\rm Sing}\,L\eq L^{[b+1]}\cup L^{[2]}$ with $L^{[b+1]}\eq{\rm Sing}\,L\cap\{xyz\eq0\}$, and the singular points of $\{xyz\eq0\}$ are disjoint from $L^{[b+1]}$, see \cite[3.4--5]{eff}. 
Assuming $b{+}1\in m\1\Z$, we have $d\eq b^2\pl 3\in m\1\Z$ if and only if $4\in m\1\Z$, that is, $m\eq 2$ or $4$. In the case $m\eq 4$, we have $b\eq 4a\mi 1$ with $a\in\Z_{>0}$.
For $a\eq1$, this arrangement is called the {\it Hessian\1} arrangement, see \cite{BDS}, \cite{Di2}, \cite{eff}, \cite{Yu}, etc. It may be called a {\it generalized Hessian\1} arrangement if $a\gess 2$. The dual $(4)$-graph has two or four connected components depending on whether $a\gess 2$ or $a\eq 1$. One connected component is always given by $\{xyz\eq0\}$, and the hypothesis of Theorem~2 is satisfied for $a\gess 2$. However, this fails in Theorem~1, since the connected component has three irreducible components.
\msn
{\bf Remark~2.6}\,(i). There are many other examples where the hypothesis of Theorem~2 is satisfied with $r'_m\gess 2$. For instance, let $L_1,L_2$ be two lines in $\PP^2$. Take $a(m{-}1)$ points
$$P_{i,j}\in L_1\stmi L_2,\q Q_{i,j}\in L_2\stmi L_1\q(i\in[1,a],\,j\in[1,m{-}1]).$$
Let $L_{i,j,j'}$ be the line passing through $P_{i,j}$, $Q_{i,j'}$ for $i\in[1,a],\,j,j'\in[1,m{-}1]$. Set
$$L=L_1\cup L_2\cup\mcup_{i,j,j'}\,L_{i,j,j'}.$$
Assume the $P_{i,j}$, $Q_{i,j'}$ are sufficiently general so that ${\rm Sing}\,L=L^{[m]}\cup L^{[2]}$ and $L^{[m]}\subset L_1\cup L_2$, where $m\gess 3$. Put $a\eq m{-}2$. Then $d=(m{-}2)(m{-}1)^2+2\in m\1\Z$, and $r'_m\eq 2$. One connected component of $L'_{\lan m\ran}$ is given by its intersection with $L_1\cup L_2$. Here $L_d$ can be $L_1$ or $L_2$, and the assumption of Theorem~1 is also satisfied in this case.
\msn
{\bf Remark~2.6}\,(ii). We have the following example where the hypothesis of Theorem~1 is satisfied with $r'_m\eq m\eq 4$. Take sufficiently general real numbers $a_i,b_i$ ($i\in[1,4]$) and a sufficiently small positive number $c$ (in particular, $a_i\pl b_j\nes a_{i'}\pl b_{j'}$ if $(i,j)\nes(i',j'))$. Consider the lines $L_{i,j,k}$ defined in $\PP^2$ by
$$y=(a_i\mi b_j\pl k\1 c)x+(a_i\pl b_j)z,$$
for $i,j\in[1,4],\,k\in\{-1,0,1\}$ with $i\nes j$, and also the lines $L_l$ ($l\eq{-1},0,1$), $L_{\infty}\,(=L_d)$ which are defined respectively by
$$x=l\1 z,\q z\eq 0.$$
They form a line arrangement $L$ of degree 40 with $r'_m\eq m\eq 4$, and the closure of three connected components of $L'_{\lan m\ran}$ with $m\eq 4$ are given by the $L_l$ ($l\eq{-1},0,1$). Note that $L$ does not support a 4-net. We may get $r'_m\eq 3$ replacing the order of $\{1,\dots,d\}$ appropriately, since $L_{\lan m\ran}$ contain irreducible connected components and $r_m\eq 4$ in the notation of the introduction with $X$ replaced by $L$, see also Remark~2.6\,(iii) just below.
\msn
{\bf Remark~2.6}\,(iii). We can modify the example in Remark~2.6\,(ii) by considering the lines $L_{i,j,k,k'}$ (instead of $L_{i,j,k}$) defined by
$$y=(a_i\mi b_j\pl k\1 c)x+(a_i\pl b_j\pl k'c')z,$$
for $i,j\in[1,4],\,k,k'\in\{-1,0,1\}$ with $i\nes j$, where $c'$ is sufficiently large compared with the $a_i,b_j$.
\sk
We have $d\eq 36\cdot 3\pl 4\eq 112$, and $r'_m\eq 4$ with $m\eq 4$. However, $r_m\eq 5$, and we may get $r'_m\eq 5$ replacing the order of $\{1,\dots,d\}$ appropriately. (In particular, this does not support a 4-net.) The assumption of Theorem~1 may be satisfied with $r'_m\eq 5$.

\end{document}